\documentclass{amsart}

\usepackage[applemac]{inputenc}   
\usepackage[T1]{fontenc}
\usepackage{lmodern}
\usepackage{graphicx}
\usepackage[english]{babel}
\usepackage{amsmath, amsfonts, amssymb, amsthm}
\usepackage{hyperref} 
\hypersetup{colorlinks = true, allcolors=blue}
\usepackage{array}
\usepackage{enumitem}

\newtheorem{thrm}{Theorem}
\newtheorem{prpstn}{Proposition}

\theoremstyle{definition}
\newtheorem{dfntn}{Definition}

\theoremstyle{remark}
\newtheorem{xmpl}{Example}
\newtheorem{rmrk}{Remark}

\begin{document}

\title{Necessary conditions for local controllability of a particular class of systems with two scalar controls} 
\author{Laetitia Giraldi}\thanks{Laetitia Giraldi and Jean-Baptiste Pomet are with McTAO team, Inria, Universit\'e C\^ote d'Azur, CNRS, LJAD, France ;\\ email: \texttt{laetitia.giraldi@inria.fr}, \texttt{jean-baptiste.pomet@inria.fr}}
\author{Pierre Lissy}\thanks{Pierre Lissy is with CEREMADE,  Universit\'e Paris-Dauphine \& CNRS UMR 7534, Universit\'e PSL, 75016 Paris, France;\\ email: \texttt{lissy@ceremade.dauphine.fr}}
\author{Cl\'ement Moreau}\thanks{Cl\'ement Moreau is with Universit\'e C\^ote d'Azur, McTAO team, Inria, CNRS, LJAD, France ; email:\texttt{clement.moreau@inria.fr}}
\author{Jean-Baptiste Pomet}
\date{\today}


\begin{abstract}
We consider affine control systems with two scalar controls, such that one control vector field vanishes at an equilibrium state. We state two necessary conditions of local controllability around this equilibrium, involving the iterated Lie brackets of the system vector fields, with controls that are either bounded, small in $\mathrm{L}^{\infty}$ or small in $\mathrm{W}^{1,\infty}$. These results are illustrated with several examples. 
\end{abstract}



\maketitle




\section{Introduction \label{section:intro}}
 
Consider a general control system of the form 
$x'(t)=f(t,x(t),u(t)),$ where $x$ is the state, $u$ is the control, $t\in [0,T]$ for some $T>0$ and $f$ is some function. Such a system is usually called \textit{controllable} if, for any two points $X$ and $Y$ in the state space, there
exists a control $u$ producing a trajectory $x$ that starts from $X$ at time $0$ and ends at $Y$ at time $T$
(see classical textbooks like \cite{Lee-Mar67,sontag2013mathematical,coron2007control}).
It is \textit{locally} controllable around a point in the state space and a value of the control
---assumed to be an equilibrium throughout all this paper, we do not discuss local controllability around a trajectory--- if
two states $X$ and $Y$ close enough to the above-mentioned equilibrium can be joined in arbitrarily small time with controls 
arbitrarily close to the reference control, by a trajectory that remains close to the equilibrium.
There are different notions of local controllability, some stronger than others, depending on the
topology used on the control, and possibly requiring that the difference with the reference control
be bounded rather than arbitrarily small, see Section~\ref{sec-defs}.

In what follows, we will restrict our attention to real analytic affine control system, defined by a finite number of real analytic vector fields.
Some sufficient conditions (see \cite{hermes1976local,hermes1982local,Suss87siam}) and some
necessary conditions (see \cite{Suss87siam,stefani1986local,Kaws88contemp,Kras98}) are given in the
literature for local controllability of these systems around an equilibrium, with a rather big gap between them that makes the subject intriguing.
These conditions all allow us to decide controllability or non-controllability based on
the value of a finite number of Lie brackets of the underlying vector fields at the equilibrium point,
i.e. on the truncation at a certain order of the series defining the real analytic vector fields; an
even more intriguing question is pointed out in \cite{Agra99}:  it is not clear whether or not, in
general (for systems lying in the above mentioned gap), a finite
number of such terms of series, or Lie brackets, is enough to decide local controllability or non controllability. See for example \cite{coron2007control} and \cite{sontag2013mathematical} for more results on the
important questions around local controllability that emerged in nonlinear control theory and for the
advances in the last decades.

\medskip

This paper is specifically concerned with control systems with two scalar inputs, of the form
\begin{equation}
  \label{eq:3}
\dot{z}=f_0 (z)  + u_1\, f_1(z)  + u_2\, f_2(z)\,,
\end{equation}
where the state $z$ is in $\mathbb{R}^n$, $f_0$, $f_1$, $f_2$
are three real analytic vector fields on $\mathbb{R}^n$ such that $f_0$
and $f_2$ vanish at the origin while $f_1$ does not:
\begin{equation}
  \label{eq:4}
  f_0(0)=0,\ f_2(0)=0,\ f_1(0)\neq0\,.
\end{equation}

Such systems have two controls but the effect of one of them vanishes at the point of interest.
In a sense, the contribution of this paper is to study to what extent the second control helps
controllability or to what extent, on the contrary, obstructions to controllability of the single
input system $\dot{z}=f_0 (z)  + u_1\, f_1(z)$ carry over when the second control  $u_2$ is
turned on.


Studying this very situation stemmed out of previous work from the authors on the controllability of
magnetic micro-swimmers \cite{giraldi2017local,giraldi2017addendum,moreau2019local}. See these
references for a description of these devices and their interest (for instance in micro-robotics and
biomedical applications). The corresponding control systems are particular cases of
\eqref{eq:3}-\eqref{eq:4}, for which
the authors have proved various controllability and non-controllability results, with the various
notions of local controllability introduced in \ref{sub:def}.

We believe that a more general treatment of systems of type \eqref{eq:3}-\eqref{eq:4}, beyond the case of magnetic micro-swimmers, is of interest to the controllability problem in control theory. It is the purpose of the present paper.


\medskip

The paper is structured as follows. Section \ref{sec-defs} is devoted to precise definitions of
various notions of local controllability and to recalling known controllability conditions for single-input systems. Our two main results are presented in Section~\ref{section:result}. Section \ref{section:examples} illustrates the results with several examples. Section \ref{section:proofs} is dedicated to the proofs. To finish, conclusions as well as some perspectives on further research are provided in Section \ref{section:conclusion}.

\section{Problem statement}
\label{sec-defs}

\subsection{Definitions of local controllability \label{sub:def}}
Let $n$ be a positive integer. Let $\mathcal{X}$ be the set of real analytic vector fields on $\mathbb{R}^n$. Let us give definitions for a general affine control system with $m$ controls:
\begin{equation}
  \label{eq:5}
\dot{z}=f_0 (z)  + \sum_{k=1}^m u_k\, f_k(z)\,.
\end{equation}
The state $z$ is in $\mathbb{R}^n$, the vector fields $f_0,\ldots,f_m$ belong to $\mathcal{X}$, the controls
$u=(u_1,\ldots,u_m)$ are in $\mathbb{R}^m$. We endow $\mathbb{R}^m$ with any norm that we denote by $|.|$. We will keep the notations $\|.\|$ for functional norms when the control is
assigned to be a function of time $t\mapsto u(t)$.

We say that $(z^{\text{eq}},u^{\text{eq}}) \in \mathbb{R}^n \times \mathbb{R}^m$ is an equilibrium point of the system if $f_0 (z^{\text{eq}}) + \sum_{k=1}^m u_k^{\text{eq}} f_k (z^{\text{eq}}) =0$.

Let us quickly review the different notions of local controllability around an equilibrium that exist in the
literature, starting with the so-called \emph{small-time local controllability (STLC)}, used by
Coron in \cite[Def. 3.2, p. 125]{coron2007control}. For $\eta \in \mathbb{R}$ such that $\eta > 0$ and $z \in \mathbb{R}^n$, we denote by $B(0,\eta)$ the open
ball for the Euclidian norm in $\mathbb{R}^n$, centered at $z$ and with radius $\eta$. 

\begin{dfntn}[STLC]
\label{def-stlc}
The control system \eqref{eq:5} is \emph{STLC at $(z^{\text{eq}},u^{\text{eq}})$} if, for every $\varepsilon >0$, there exists $\eta >0$ such that, for every $z_0,z_1$ in $B(z^{\text{eq}},\eta)$, there exists a control $u(\cdot)$ in $ \mathrm{L}^{\infty}([0,\varepsilon],\mathbb{R}^m)$ such that the solution of the control system $z(\cdot):[0,\varepsilon]\to\mathbb{R}^{n}$ of \eqref{eq:5} 
satisfies $z(0)=z_0$, $z(\varepsilon)=z_1$, and
\begin{equation*}
\| u - u^{\text{eq}} \| _{\mathrm{L}^{\infty}([0,\varepsilon],\mathbb{R}^m)} \leqslant \varepsilon\,.
\end{equation*}
\end{dfntn}

Note that this notion requires the time to be arbitrarily small and the control to be arbitrarily close to the equilibrium control. Nevertheless, another notion that only requires boundedness of the control can be found in the works of Hermes \cite{hermes1982local} and Sussmann \cite{sussmann1983lie}\footnote{The exact definition given in \cite{sussmann1983lie} supposes an \textit{a priori} bound on the control, uses the notion of reachable space, and is hence written in a more condensed manner. We rephrase it here to match the structure of the first definition.} among others. This second notion, while not equivalent to the first one, is sometimes called STLC as well. In order to avoid the confusion, we will call it $\alpha$-STLC according to the following definition:

\begin{dfntn}[$\alpha$-STLC]
\label{def-stlcq}
Let $\alpha \geqslant 0$.
The control system \eqref{eq:5} is \emph{$\alpha$-STLC at $(z^{\text{eq}},u^{\text{eq}})$} 
if, for every $\varepsilon >0$, there exists $\eta >0$ such that, for every $z_0,z_1$ in $B(z^{\text{eq}},\eta)$,
there exists a control $u(\cdot)$ in $\mathrm{L}^{\infty}([0,\varepsilon],\mathbb{R}^m)$ such that the solution of the control system $z(\cdot):[0,\varepsilon]\to\mathbb{R}^{n}$ of \eqref{eq:5} 
satisfies $z(0)=z_0$, $z(\varepsilon)=z_1$, and
\begin{equation*}
\| u -u^{\text{eq}} \| _{\mathrm{L}^{\infty}([0,\varepsilon],\mathbb{R}^m)} \leqslant \alpha + \varepsilon\,.
\end{equation*}
\end{dfntn}

\begin{rmrk}  We can easily see that 0-STLC is then identical to STLC. If $\alpha >0$, the second notion is weaker than the first one, as the norm of the control can remain ``far'' from the equilibrium control as the ball radius $\eta$ gets arbitrary small. 
\label{rmk:2} \end{rmrk}

\begin{rmrk} For a given control system, the smallest possible value of $\alpha$ depends on the norm $| \cdot |$ chosen for the control. However, it does not depend on the norm we put on the state space, justifying the  choice of a particular norm on $\mathbb R^n$.
\end{rmrk} 

More recently, a new notion has been introduced by Beauchard and Marbach in \cite[Definition 4]{beauchard2018quadratic}. The idea is to ensure the smallness, not only of the control, but also of its derivatives. Hence its norm will be bounded in the Sobolev spaces $\mathrm{W}^{k,\infty}$.

\begin{dfntn}[$\mathrm{W}^{k,\infty}$-STLC]
\label{def-stlcw}
Let $k \in \mathbb{N}$.
The control system \eqref{eq:5} is \emph{$\mathrm{W}^{k,\infty}$-STLC at $(z^{\text{eq}},u^{\text{eq}})$}
if, for every $\varepsilon >0$, there exists $\eta >0$ such that, for every $z_0,z_1$ in $B(z^{\text{eq}},\eta)$,
there exists a control $u(\cdot )$ in $\mathrm{W}^{k,\infty}([0,\varepsilon],\mathbb{R}^m)$ such that the solution of the control system $z(\cdot):[0,\varepsilon]\to\mathbb{R}^{n}$ of \eqref{eq:5} 
satisfies $z(0)=z_0$, $z(\varepsilon)=z_1$, and
\begin{equation*}
\| u -u^{\text{eq}} \| _{\mathrm{W}^{k,\infty}([0,\varepsilon],\mathbb{R}^m)} \leqslant \varepsilon\,.
\end{equation*}
\end{dfntn}

\begin{rmrk} Like $0$-STLC in Remark \ref{rmk:2}, $\mathrm{W}^{0,\infty}$-STLC is identical to STLC. When $k>0$, $\mathrm{W}^{k,\infty}$-STLC is stronger than STLC, because it requires the control to be sufficiently smooth. \end{rmrk}


\begin{rmrk} The different STLC notions can be ordered in an implication chain. With $0 < \alpha_1 < \alpha_2$ and $k_1 < k_2$ in $\mathbb{N}$, and for a given norm on the control space, one has
\[
\mathrm{W}^{k_2,\infty}\text{-STLC} \Rightarrow \mathrm{W}^{k_1,\infty}\text{-STLC} \Rightarrow \text{STLC} \Rightarrow \alpha_1\text{-STLC} \Rightarrow \alpha_2\text{-STLC}
\]
\end{rmrk}

\subsection{Known results for single-input systems}

Consider an affine control system like \eqref{eq:5} with $m=1$:
\begin{equation} 
\dot{z}=f_0 (z) + u_1(t) f_1(z)
\label{eq:1}
\end{equation}
with $z$ in $\mathbb{R}^n$, $f_0, f_1$ in $\mathcal{X}$ and $u_1$ a control function in $\mathrm{L}^1([0,T])$. For some $T_u$ in $]0,T]$, \eqref{eq:1} admits a unique maximal solution (see e.g. \cite[Proposition 2]{beauchard2018quadratic}). 
Up to a translation, we can assume that $(0,0)$ is an equilibrium of \eqref{eq:1} (which means in particular that $f_0 (0) = 0$). 

Since $m=1$, we can assume here without loss of generality that $|\cdot|$ is the usual absolute value.

If $f$ and $g$ are two vector fields in $\mathcal{X}$, $[f,g]$ denotes the Lie bracket of $f$ and $g$ and $\mathrm{ad}_{f}^k g$ is defined by induction with $\mathrm{ad}_{f}^0g=g$ and $\mathrm{ad}_{f}^k g=[f,\mathrm{ad}_{f}^{k-1}g]$.

For $k \in \mathbb{N}$, we denote by $\mathrm{Lie}(f_0,f_1)$ the Lie algebra of vector fields
generated by $f_0$ and $f_1$, $S_k$ the subspace of $\mathcal{X}$ spanned by all the Lie
brackets of $f_0, f_1$ containing $f_1$ at most $k$ times, and $S_k (0)$ the subspace of $\mathbb{R}^n$
spanned by the value at 0 of the elements of $S_k$.

A sufficient condition for STLC is given by the following proposition.

\begin{prpstn}[\cite{sussmann1983lie}, Theorem 2.1, p.688]
\label{sufficient-condition} If 
\begin{equation}
\label{eq:lie-cond}
\{ g(0), g \in \mathrm{Lie}(f_0,f_1)\} = \mathbb{R}^n,
\end{equation}
and, for all $k$ in $\mathbb{N}$, 
\begin{equation}
\label{eq:lie-cond-2}
S_{2k+2}(0) \subset S_{2k+1}(0),
\end{equation} then system \eqref{eq:1} is STLC.
\end{prpstn}

A partial converse result states that the condition \eqref{eq:lie-cond} is necessary for any form of $STLC$. 

Condition \eqref{eq:lie-cond-2} is violated if for some $k \in \mathbb{N}$, some brackets in $S_{2k+2}$, once evaluated at $0$, do not belong to $S_{2k+1}(0)$. We call such brackets \textit{bad} ones, since they constitute a potential \textit{obstruction} to local controllability. The bad brackets can be seen as directions towards which the system drifts, thus potentially (but not necessarily) preventing local controllability. The fact that they do not belong to $S_{2k+1}(0)$ at $0$ for some $k$ means that they do not share directions with the brackets in $S_{2k+1}$; hence those brackets evaluated at 0 cannot be used to ``compensate'' the drift induced by the bad brackets. Dealing with this issue can be done by finding other brackets sharing directions at $0$ with the bad brackets, in order to compensate the drift. We say that those new brackets help to \textit{neutralize} the bad ones. 

The lowest-order possible obstruction occurs if the bad bracket in $S_2$
\begin{equation}
B_1 = [f_1,[f_0,f_1]]
\label{eq:first-bad}
\end{equation}
is such that $B_1(0)$ does not belong to $S_1(0)$. We call $B_1$ the \textit{first bad bracket}.
The following result has been shown by Sussmann in \cite[Proposition 6.3, p.707]{sussmann1983lie}:

\begin{prpstn}
\label{thm:sussmann}
Assume $f_0(0) = 0$ and $B_1(0)  \not \in S_1 (0)$. Then,  for any $\alpha \geqslant  0$, \eqref{eq:1} is not $\alpha$-STLC.
\end{prpstn}

If $B_1(0)$ does belong to $S_1(0)$, the next lowest-order bracket in $S_2$ that can obstruct controllability is what we call the \textit{second bad bracket} and denote by $B_2$: 
\begin{equation}
B_2=[[f_0,f_1],[f_0,[f_0,f_1]]].
\label{eq:second-bad}
\end{equation}

For scalar-input systems, Sussmann noticed in \cite[p.710]{sussmann1983lie} that one may or may not get STLC under the hypothesis $B_2(0) \not \in S_1(0)$. In \cite{Kaws88contemp}, Kawski obtained a new necessary condition by refining the space $S_1$:

\begin{prpstn}
Let $S' = \mathrm{Span} ( \{ \mathrm{ad}^{k}_{f_0} ( \mathrm{ad}^3_{f_1} f_0 ), k \in \mathbb{N} \} )$.

If $B_2(0) \not \in S_1 (0) + S'(0)$, then,  for any $\alpha \geqslant  0$, \eqref{eq:1} is not $\alpha$-STLC.
\label{thm:kawski}
\end{prpstn}

More recently, in \cite[Theorem 3]{beauchard2018quadratic}, Beauchard and Marbach showed another result by using another notion of local controllability.

\begin{prpstn}
If $B_2(0) \not \in S_1(0)$, then \eqref{eq:1} is not $\mathrm{W}^{1,\infty}$-STLC. 
\label{thm:beauchard}
\end{prpstn}

Proposition \ref{thm:kawski} states that $B_2$ can be neutralized if it shares its direction at $0$ with a particular class of brackets in $S_3$, while Proposition \ref{thm:beauchard} states that $B_2$ can only be neutralized if the derivative of the control is ``not too small''.

Concerning systems with control in $\mathbb{R}^m, m \geqslant 2$, a general sufficient condition for local controllability, in the vein of Proposition \ref{sufficient-condition} but more complex, can be found in \cite{sussmann1987general}, but no necessary condition is known, to the best of our knowledge.  The main results of this paper, stated in the next section, are a step in this direction in that they give an extension of the necessary conditions contained in Propositions \ref{thm:sussmann}, \ref{thm:kawski} and \ref{thm:beauchard} to  the case where the system has two scalar controls, and the vector field associated to the second control vanishes.


\section{Main results \label{section:result}}

We now consider the affine control system \eqref{eq:3} (which is also system \eqref{eq:5} with $m=2$):
\begin{equation*}
\dot{z}=f_0 (z) + u_1(t) f_1(z) + u_2 (t) f_2 (z),
\end{equation*}
with $z \in \mathbb{R}^n$, $f_0$, $f_1$, $f_2$ in $\mathcal{X}$ and $u_1$, $u_2$ control functions in $\mathrm{L}^1([0,T])$. 

We assume that \eqref{eq:4} is verified, i.e. $f_0 (0) = 0, f_2(0) = 0, f_1(0) \neq 0$, and we study local controllability for $(z,(u_1,u_2))$ close to the equilibria $(0,(0,u_2^{\text{eq}}))$, with $u_2^{\text{eq}}$ arbitrary.

\begin{rmrk}
  A more general situation than \eqref{eq:4} would be to consider an equilibrium $(z^{\text{eq}},(u_1^{\text{eq}},u_2^{\text{eq}}))$ (i.e. $f_0 (z^{\text{eq}})  + u_1^{\text{eq}}\,
  f_1(z^{\text{eq}})  + u_2^{\text{eq}}\, f_2(z^{\text{eq}})=0$) such that the rank of
  $\{f_1(z^{\text{eq}}), f_2(z^{\text{eq}})\}$ is 1. In that case, one may recover \eqref{eq:4} by defining new variables and controls $(Z,U_1,U_2)$ with the linear transformation $z=Z+z^{\text{eq}},
  u_1=\lambda_2 U_1+\lambda_1 U_2, u_2=-\lambda_1 U_1+\lambda_2 U_2$ where $(\lambda_1,\lambda_2)$
  nonzero such that $\lambda_1 \, f_1(z^{\text{eq}})  +\lambda_2\, f_2(z^{\text{eq}})=0$. This transformation brings us back to the study of a system of type \eqref{eq:3}-\eqref{eq:4}.
\end{rmrk}

Let $R_1$ be the subspace of $\mathcal{X}$ spanned by all the iterated Lie brackets of $f_0, f_1, f_2$ containing $f_1$ at most one time, and $R_1 (0)$ the subspace of $\mathbb{R}^n$ spanned by the value at 0 of the elements of $R_1$. 

\subsection{Obstruction coming from the first bad bracket}

Let us now state our main results. Recall that $B_1=[f_1,[f_0,f_1]]$ and $B_2=[[f_0,f_1],[f_0,[f_0,f_1]]]$.



\begin{thrm}
  \label{main-thm-bis}
  Consider system \eqref{eq:3} under Assumption \eqref{eq:4}. Assume $B_1 (0) \not \in R_1 (0)$.
\begin{enumerate}[label=\arabic*.]
\item \label{first-case} If $B_1(0) \in R_1(0)+\mathrm{Span}([f_1, [f_2, f_1]](0))$, let $\beta\in\mathbb R$ be such that
  \begin{equation*}
    B_1(0) + \beta [f_1,[f_2,f_1]](0)\in R_1(0).
  \end{equation*}
Then, for any $u_2^{\text{eq}}\in\mathbb R$ such
  that $u_2^{\text{eq}}\neq\beta$,
 system \eqref{eq:3} is not STLC at $(0,(0,u_2^{\text{eq}}))$.
  \item \label{second-case} If $B_1(0)\not\in R_1(0)+\mathrm{Span}([f_1,[f_2,f_1]](0))$, then, for any $u_2^{\text{eq}}\in\mathbb R$
  and any $\alpha\geqslant 0$,
  system \eqref{eq:3}  is not $\alpha$-STLC at $(0,(0,u_2^{\text{eq}}))$. \end{enumerate}
\end{thrm}

\begin{rmrk} In case \ref{second-case}, the second control does not improve controllability with respect to the single-input system obtained by taking $u_2 = 0$.
\end{rmrk}

\begin{rmrk} In case \ref{first-case}, the fact that the brackets $[f_1, [f_2, f_1]](0)$ and $B_1(0)$ share a common direction is crucial. It allows the bracket $[f_1, [f_2, f_1]]$ to possibly neutralize the bad bracket $B_1$ through the particular control $u_2^{\text{eq}}=\beta$. This critical value of the control is the only value around which system \eqref{eq:3} may be STLC. \end{rmrk}

\begin{rmrk}
In \cite{giraldi2017addendum}, a result similar to case \ref{first-case} is shown for a particular system of type \eqref{eq:3}-\eqref{eq:4} describing the movement of a magnetized micro-swimmer. This particular result led to the generalizations presented in this paper. 
Furthermore, the proof of Theorem \ref{main-thm-bis} is based on the existence of a suitable local change of coordinates, that is performed explicitly in \cite{giraldi2017addendum} for the micro-swimmer system.
\end{rmrk}

\begin{rmrk} Up to a translation, one can always study controllability around the null equilibrium $(0,(0,0))$. Let us define the affine feedback transformation on the control $u_2$: $\tilde{u}_2=u_2 - \beta$. With this transformed control, system \eqref{eq:3} becomes
\begin{equation}
\dot{z} = \tilde{f}_0 (z) + u_1 \tilde{f}_1 (z) + \tilde{u}_2 \tilde{f}_2 (z)
\label{eq:2bis}
\end{equation}
with $\tilde{f}_0 = f_0 + \beta f_2$, $\tilde{f}_1 = f_1$ and $\tilde{f}_2 = f_2$.

Note that $[\tilde{f}_1,[\tilde{f}_0,\tilde{f}_1]](0)=[f_1,[f_0,f_1]](0)+\beta [f_1,[f_2,f_1]](0) \in R_1(0)$. Assume that system \eqref{eq:2bis} is STLC at $(0,(0,0))$. Let $\varepsilon$ be a positive real number. Let $\eta$ be the associated parameter from Definition \ref{def-stlc}, and $z_0$, $z_1$ in $B(0,\eta)$. There exists controls $u_1$ and $\tilde{u}_2$ in $\mathrm{L}^{\infty}([0,\varepsilon])$ such that the solution of \eqref{eq:2bis} with $z(0)=z_0$ and these controls verify $z(\varepsilon)=z_1$, and $$\| u_1 \| _{\mathrm{L}^{\infty}([0,\varepsilon],\mathbb{R})} \leqslant \varepsilon\,, \| \tilde{u}_2 \| _{\mathrm{L}^{\infty}([0,\varepsilon],\mathbb{R})} \leqslant \varepsilon.$$ Hence, the solution of system \eqref{eq:3} with $z(0)=z_0$ and controls $u_2 = \beta + \tilde{u}_2$ and $u_1$ verifies $z(\varepsilon)=z_1$. Moreover, $\| u_2 - \beta \| _{\mathrm{L}^{\infty}([0,\varepsilon],\mathbb{R})} \leqslant  \varepsilon$ and $\| u_1 \| _{\mathrm{L}^{\infty}([0,\varepsilon],\mathbb{R})} \leqslant \varepsilon$.

Therefore, if system \eqref{eq:2bis} is STLC at $(0,(0,0))$, then system \eqref{eq:3} is STLC at $(0,(0,\beta))$.
\label{remark:S'}
\end{rmrk}

\subsection{Obstruction coming from the second bad bracket}

In order to state our result, let us introduce a complementary notion of local controllability, fit with the type of systems \eqref{eq:3}-\eqref{eq:4} we are interested in.

\begin{dfntn}
\label{def-stlc2}
Let $k \in \mathbb{N}$ and $\alpha$ in $\mathbb{R}$ such that $\alpha \geqslant 0$.
The control system \eqref{eq:3} is $(\mathrm{W}^{k,\infty},\alpha)$-STLC at $(z^{\text{eq}},(u_1^{\text{eq}},u_2^{\text{eq}}))$
if, for every $\varepsilon >0$, there exists $\eta >0$ such that, for every $z_0,z_1$ in $B(z^{\text{eq}},\eta)$,
there exists a control $(u_1(\cdot ),u_2(\cdot))$ in $\mathrm{W}^{1,\infty}([0,\varepsilon],\mathbb{R}) \times \mathrm{L}^{\infty}([0,\varepsilon],\mathbb{R})$ such that the solution of the control system $z(\cdot):[0,\varepsilon] \to \mathbb{R}^{n}$ of \eqref{eq:3} 
satisfies $z(0)=z_0$, $z(\varepsilon)=z_1$, and
\begin{equation*}
\| u_1 -u_1^{\text{eq}} \| _{\mathrm{W}^{k,\infty}([0,\varepsilon],\mathbb{R})} \leqslant \varepsilon, \quad \| u_2 -u_2^{\text{eq}} \| _{\mathrm{L}^{\infty}([0,\varepsilon],\mathbb{R})} \leqslant \alpha + \varepsilon.
\end{equation*}
\end{dfntn}

\begin{rmrk}
The norms used for each control are different in this STLC notion. It fits the nature of system \eqref{eq:3}, where the second control plays a particular role due to the fact that $f_2$ vanishes at $0$. This could be seen as a form of ``hybrid'' small-time local controllability.
\end{rmrk}


\begin{thrm}
Consider system \eqref{eq:3} under assumption \eqref{eq:4}. Assume that $B_1(0)\in R_1(0)$ and that $B_2(0)\not\in\mathrm{Span}(R_1(0), [f_1,[f_2,f_1]](0))$. 
\begin{enumerate}[label=\arabic*.]
\item \label{prop-first-case} If $B_2(0) \in \mathrm{Span} ( R_1(0), \{ [[f_i,f_1],[f_j,[f_k,f_1]]](0), (i,j,k) \in \{ 0,2 \}^3, (i,j,k) \neq (0,0,0) \} )$, then system \eqref{eq:3} is not $(\mathrm{W}^{1,\infty},0)$-STLC at $(0,(0,0))$.
\item \label{prop-second-case} Else,  for any $\alpha \geqslant  0$, system \eqref{eq:3} is not $(\mathrm{W}^{1,\infty},\alpha)$-STLC at $(0,(0,0))$.
\end{enumerate}
\label{thm:b2}
\end{thrm}

\begin{rmrk} The natural hypothesis, instead of $B_2(0) \not \in \mathrm{Span}(R_1(0), [f_1,[f_2,f_1]](0))$, would be $B_2(0) \not \in R_1(0)$. It will become clear in the proof why we need to strengthen this hypothesis. The case $B_2(0) \in \mathrm{Span}(R_1(0), [f_1,[f_2,f_1]](0))$ is still under our investigation; see Example \ref{example:2-3} for more details.
\end{rmrk}

\section{Illustrating examples and applications \label{section:examples}}

\subsection{Examples for the first bracket obstruction}

\subsubsection{Case where the second control cannot help to neutralize $B_1$}

In case \ref{second-case} of Theorem \ref{main-thm-bis}, the second control $u_2$ cannot neutralize the obstruction to local controllability induced by $B_1$. The following example illustrates that case.

\begin{xmpl} Consider the system
\begin{equation} 
\left \{
\begin{array}{r l}
\dot{x} & \displaystyle = y^2 + y u_1, \\ 
\dot{y} & \displaystyle = 2y - u_1 + x u_2.
\end{array}
\right . 
\label{example-0}
\end{equation}
It is of the form \eqref{eq:3} with
\[
f_0 = \begin{pmatrix} y^2 \\ 2 y \end{pmatrix}, \quad f_1 = \begin{pmatrix} y \\ -1 \end{pmatrix}, \quad f_2 = \begin{pmatrix} 0 \\ x  \end{pmatrix}.
\]
Straightforward computations show that
\[
R_1 (0) = \mathrm{Span}(\mathbf{e}_2), \quad [f_1,[f_0,f_1]](0)=-6 \mathbf{e}_1, \quad [f_1,[f_2,f_1]](0)= \mathbf{e}_2,
\]
so we are in the case \ref{second-case} of Theorem \ref{main-thm-bis}. Therefore,  for any $\alpha \geqslant 0$ and any $u_2^{\text{eq}} \in \mathbb{R}$, system \eqref{example-0} is not $\alpha$-STLC at $(0,(0,u_2^{\text{eq}}))$.
\end{xmpl}

\subsubsection{Case where STLC is retrieved thanks to the second control}

In case \ref{first-case}, Theorem \ref{main-thm-bis} states that the system is not STLC around the equilibria $(0,(0,u_2^{\text{eq}}))$, unless $u_{\text{eq}}$ is equal to a particular value $\beta$, that allows the bracket $[f_1,[f_2,f_1]]$ to neutralize the bad bracket $B_1$. Around the equilibrium $(0,(0,\beta))$, the system can then be STLC, like in the next example. The method used in the following example to show STLC was introduced in \cite{moreau2019local} to show local controllability of magnetically driven micro-swimming robots. We reproduce it here on a simpler system.

\begin{xmpl}  \label{example:swimmer} Consider the system
\begin{equation} 
\left \{
\begin{array}{r l}
\dot{x} & \displaystyle = y^2 + y u_1 - \frac{2}{\alpha} y^2 u_2,  \vspace{2mm} \\ 
\dot{y} & \displaystyle = 2y - u_1 - \frac{1}{\alpha} y u_2,
\end{array}
\right . 
\label{example-1}
\end{equation}
for some  $\alpha \neq 0$. Here we have
\[
f_0 = \begin{pmatrix} y^2 \\ 2 y \end{pmatrix}, \quad f_1 = \begin{pmatrix} y \\ -1 \end{pmatrix}, \quad f_2 = -\frac{1}{\alpha} \begin{pmatrix} 2 y^2 \\ y  \end{pmatrix}.
\]
Straightforward computations show that
\begin{equation*}
R_1 (0) = \mathrm{Span}(\mathbf{e}_2), \quad [f_1,[f_0,f_1]](0)=-6 \mathbf{e}_1, \quad [f_1,[f_2,f_1]](0)=\frac{6}{\alpha} \mathbf{e}_1,
\end{equation*}
so we are in the case \ref{first-case} of Theorem \ref{main-thm-bis}. Therefore, the system \eqref{example-1} is not STLC at $(0,(0,u_2^{\text{eq}})$ for any $u_2^{\text{eq}} \neq \alpha$. 

As in Remark \ref{remark:S'}, we define the feedback control $\tilde{u}_2$ such that $u_2 = \alpha + \tilde{u}_2$, that neutralizes the bracket $B_1$. With this control, the transformed system \ref{example-1} reads
\begin{equation} 
\left \{
\begin{array}{r l}
\dot{x} & \displaystyle = - y^2  + y u_1 - \frac{2}{\alpha} y^2 \tilde{u}_2, \vspace{2mm} \\ 
\dot{y} & \displaystyle = y  - u_1- \frac{1}{\alpha} y \tilde{u}_2,
\end{array}
\right . \quad \text{so that } \tilde{f}_0 = \begin{pmatrix} -y^2 \\ y \end{pmatrix}.
\label{example-1-bis}
\end{equation}
Let us show that this system is STLC at $(0,(0,0))$. To this end, we use the sufficient Sussmann condition for controllability \cite[Theorem 7.3]{sussmann1987general} with $\theta=1$ and the notation for $G_{\eta}$ introduced in \cite[Definition III.10]{giraldi2017local}. Since $[f_1,[f_2,f_1]](0)=\frac{6}{\alpha} \mathbf{e}_1$, the Lie brackets of order $3$ generate the whole space, i.e. $G_{\eta}$ is the whole tangent space if $\eta > 3$. The only Lie brackets of order at most $3$ with an even number of $1$ and $2$ are $[f_1,[\tilde{f}_0,f_1]]$ and $[f_2,[\tilde{f}_0,f_2]]$, which are both zero and therefore belong trivially to $G_3$. 

Hence, the Sussmann condition from \cite{sussmann1987general} is verified and system \eqref{example-1-bis} is STLC at $(0,(0,0))$. We conclude that the system \eqref{example-1} is STLC at $(0,(0,\alpha))$ (see Remark \ref{remark:S'} for details).
\end{xmpl} 

\subsubsection{Application to micro-swimmer robots \label{sssection:microswimmer}}

The present paper was motivated by the work on controllability of micro-swimmer robot models made in \cite{giraldi2017local,giraldi2017addendum,moreau2019local}. The two swimmers studied in these papers are made of two (respectively three) magnetized rigid segments, linked together with torsional springs, immersed in a low-Reynolds number fluid, and driven by a uniform in space, time-varying magnetic field $\mathbf{H}$. The swimmers' movement is assumed to be planar. The magnetic field $\mathbf{H}$ belongs to the swimmers' plane and can therefore be decomposed, in the moving basis associated to the first segment, in two components called $(H_{\bot},H_{\parallel})$.

Seeing the magnetic field as a control function, the dynamics of both swimmers write as control systems that are exactly of type \eqref{eq:3}-\eqref{eq:4}:
\begin{equation}
\dot{\mathbf{z}} = f_0 (\mathbf{z}) + H_{\bot} f_1(\mathbf{z}) + H_{\parallel} f_2(\mathbf{z}),
\label{eq:ex-swimmer}
\end{equation}
with the state $\mathbf{z}$ in $\mathbb{R}^4$ for the two-link swimmer (resp. $\mathbb{R}^5$ for the three-link swimmer). The detailed expressions of $f_0$, $f_1$ and $f_2$ with respect to the system parameters are given in \cite[Equations (12) to (16)]{giraldi2017local} (resp. \cite[Appendix]{moreau2019local}).

Moreover, assumptions \eqref{eq:4} are verified. Hence, for all $H_{\parallel}$ in $\mathbb{R}$, $(0,(0,H_{\parallel}))$ is an equilibrium point (the first zero is short for $(0,0,0,0)$ in $\mathbb{R}^4$ (resp. $(0,0,0,0,0)$ in $\mathbb{R}^5$)). One also has $R_1(0) = \mathrm{Span}(\mathbf{e}_2,\mathbf{e}_3,\mathbf{e}_4)$ (resp. $R_1(0) = \mathrm{Span}(\mathbf{e}_2,\mathbf{e}_3,\mathbf{e}_4,\mathbf{e}_5)$) and the brackets of interest for Theorem \ref{main-thm-bis} read:
\begin{equation*}
[f_1,[f_0,f_1]](0) = (a_2,0,0,0) \quad \text{(resp. } [f_1,[f_0,f_1]](0) = (a_3,0,0,0,0) \text{)}
\end{equation*}
and
\begin{equation*}
[f_1,[f_2,f_1]](0) = (b_2,0,0,0) \quad \text{(resp. } [f_1,[f_2,f_1]](0) = (b_3,0,0,0,0) \text{)}
\end{equation*}
with $a_2, a_3, b_2, b_3$ constants that are nonzero under generic assumptions on the system parameters -- see \cite[Assumption III.2]{giraldi2017local} (resp. \cite[Assumption 1]{moreau2019local}).

We can therefore apply Theorem \ref{main-thm-bis}, case \ref{first-case} and conclude that the two-link swimmer (resp. three-link swimmer) is not STLC at $(0,(0,H_{\parallel}))$ for any $H_{\parallel}$ such that $H_{\parallel} \neq \frac{a_2}{b_2} $ (resp. $H_{\parallel} \neq \frac{a_3}{b_3}$). 

In \cite{moreau2019local}, it is shown that the two-link swimmer (resp. the three-link swimmer) is indeed STLC at $(0,(0,\frac{a_2}{b_2}))$ (resp. $(0,(0,\frac{a_3}{b_3})))$, using the technique displayed in Example \ref{example:swimmer}. However, the question of STLC at \emph{other} equilibria of type $(0,(0,H_{\parallel}))$ was left open in \cite[Remark 5]{moreau2019local}. Theorem \ref{main-thm-bis} allows to answer that question: $(0,(0,\frac{a_2}{b_2}))$ (resp. $(0,(0,\frac{a_3}{b_3})))$ is the \emph{only} equilibrium of this type for which the swimmer is STLC. 

\begin{rmrk}
Former studies on the two-link swimmer had led to the following results: in \cite{giraldi2017local}, it is shown that the control system \eqref{eq:ex-swimmer} associated to the 2-link swimmer is $\left ( 2 \frac{a_2}{b_2} \right )$-STLC at $(0,(0,0))$; in \cite{giraldi2017addendum}, it is shown that it is moreover not STLC at $(0,(0,0))$. The proof of this last result features an explicit construction of the function $\Phi$ that is used in the proof of Theorem \ref{main-thm-bis} below. 
\end{rmrk}

\subsection{Examples for the second bracket obstruction}

We start by recalling the classical scalar-input example given by Sussmann in \cite[Equation (6.12), p. 711]{sussmann1983lie}:
\begin{equation*} 
\left \{
\begin{array}{r l}
\dot{x} & = u_1, \\
\dot{y} & = x, \\
\dot{z} & = x^3 + y^2.
\end{array}
\right.
\end{equation*}
This system is STLC in spite of $B_2(0)$ being outside of $S_1(0)$. Proposition \ref{thm:beauchard} has shown that it is nonetheless not $\mathrm{W}^{1,\infty}$-STLC: hence it can only be controlled around $(0,0)$ with small controls if those control's derivatives are ``not too small''. 

In the following examples, we add a second control to this system and look at the effect of this second control on controllability, i.e. if the second control can or cannot help neutralize the bad bracket $B_2$ to retrieve controllability with $u_1$ small in $\mathrm{W}^{1,\infty}$. 

From now on, we consider the control system
\begin{equation}  \label{example-2-1}
\left \{
\begin{array}{r l}
\dot{x} & = u_1, \\
\dot{y} & = x + \phi (x,y,z) u_2, \\
\dot{z} & = x^3 + y^2 + \psi(x,y,z) u_2,
\end{array}
\right.
\end{equation}
with $\phi$ and $\psi$ real analytic functions from $\mathbb{R}^3$ to $\mathbb{R}$ that vanish at $(0,0,0)$.

Straightforward computations show that $B_1(0) = 0$ (so it trivially belongs to $R_1(0)$), and $B_2(0) = - 2 \mathbf{e}_3$.

\subsubsection{Existence of several neutralizing brackets}

Case \ref{prop-second-case} of Theorem \ref{thm:b2} states that the second control cannot improve controllability, even if it is not small. On the other hand, in case \ref{prop-first-case}, some brackets can neutralize the bad bracket $B_2$. In Theorem \ref{main-thm-bis}, only the bracket $[f_1,[f_2,f_1]]$ can play this role, whereas several brackets can in Theorem \ref{thm:b2}. This is illustrated in the following example, where we choose different sets of functions $\phi$ and $\psi$ in order to make the bad bracket $B_2$ colinear to several different brackets at 0. 

\begin{xmpl}\label{example:2-2} 
Take $\phi (x,y,z) = 0$ and $\psi (x,y,z) = y^2$. Then
\[
f_2 = \begin{pmatrix} 0 \\ 0 \\ y^2 \end{pmatrix},
\]
and one can compute that $R_1(0) = \mathrm{Span}(\mathbf{e}_1,\mathbf{e}_2)$, $[f_1,[f_2,f_1]](0)=0$ and $[[f_0,f_1],[f_2,[f_0,f_1]]](0)=- 2 \mathbf{e}_3$. We are in the case \ref{prop-first-case} of Theorem \ref{thm:b2}: the system \eqref{example-2-1} is  not $(\mathrm{W}^{1,\infty},0)$-STLC in the sense of Definition \ref{def-stlc2}.

Through smart choices of $\phi$ and $\psi$, one can make any other bracket of type $[[f_i,f_1],[f_j,[f_k,f_1]]]$ be the one that is colinear to $B_2$. For instance, choosing $\phi (x,y,z) = x$ and $\psi (x,y,z) = xy$ gives $R_1(0) = \mathrm{Span}(\mathbf{e}_1,\mathbf{e}_2)$, $[f_1,[f_2,f_1]](0)=0$ and $[[f_2,f_1],[f_0,[f_2,f_1]]](0)= - 2 \mathbf{e}_3$ for system \eqref{example-2-1}.
\end{xmpl} 

\subsubsection{Role of the bracket $[f_1,[f_2,f_1]]$}

Theorem \ref{thm:b2} requires that $B_2(0) \not \in \mathrm{Span}(R_1(0), [f_1,[f_2,f_1]](0))$. When this hypothesis is not verified, then the bracket $[f_1,[f_2,f_1]]$ may be used as well to neutralize $B_2$ and retrieve $\mathrm{W}^{1,\infty}$-STLC. The following example is an illustration of this fact.

\begin{xmpl}\label{example:2-3} Consider the control system \eqref{example-2-1} with $\phi (x,y,z)=0$ and $\psi(x,y,z) = x^2$.
One can check that: $R_1(0) = \mathrm{Span} (\mathbf{e}_1,\mathbf{e}_2)$, $B_1(0) = 0$ (so it trivially belongs to $R_1(0)$), and $B_2(0) = 2 \mathbf{e}_3$. Moreover, $[f_1,[f_2,f_1]](0) = 2 \mathbf{e}_3$, so $B_2(0) \in \mathrm{Span}(R_1(0), [f_1,[f_2,f_1]](0))$.

Let us show that in this case, system \eqref{example-2-1} is in fact $\mathrm{W}^{1,\infty}$-STLC. Let $\varepsilon$ be a real positive number. The first step, inspired by the return method of Coron \cite[Chapter 6]{coron2007control}, is to construct a loop trajectory, that goes from and back to $0$. We define $u_1$ and $u_2$ on $[0, 2 \pi \varepsilon]$ by 
\begin{equation}
\label{eq:u1u2}
\begin{array}{r l}
u_1 (t) & = \frac{1}{4} \varepsilon^2 \sin \left ( \frac{t}{\varepsilon} \right ) - \frac{1}{2} \varepsilon^2 \sin \left ( \frac{2t}{\varepsilon} \right ); \\
u_2 (t) & = \frac{5}{16} \varepsilon^2 \cos \left ( \frac{t}{\varepsilon} \right ) - \frac{3}{4} \varepsilon^2 \cos \left ( \frac{2t}{\varepsilon} \right ).
\end{array}
\end{equation}
One can check that, with these controls, the solution of \eqref{example-2-1} starting at $(0,0,0)$ verifies $x(2 \pi \varepsilon) = y(2 \pi \varepsilon) = z(2 \pi \varepsilon) = 0$. Now let us show that a small perturbation of this loop trajectory allows to access a neighbourhood of 0. We define the perturbed control 
\begin{equation}
\label{eq:u1per}
\textstyle u_{1,\mathrm{per}} (t) = \frac{1}{4} \varepsilon^2 \sin \left ( \frac{t}{\varepsilon} \right ) - \frac{1}{2} \varepsilon^2 \sin \left ( \frac{2t}{\varepsilon} \right ) + a + bt + ct^3
\end{equation}
for $(a,b,c) \in \mathbb{R}^3$ such that $|a|, |b|$ and $|c|$ are smaller than $\varepsilon$. Note that for all $|a|, |b|$ and $|c|$ small enough, $\| u_{1,\mathrm{per}} \|_{\mathrm{W}^{1,\infty}} \leqslant \varepsilon$ and $\| u_2 \| _{\mathrm{W}^{1,\infty}} \leqslant \varepsilon$. Let $F : \mathbb{R}^3 \to \mathbb{R}^3$ be the application that maps $(a,b,c)$ to \newline $(x_{\mathrm{per}}(2 \pi \varepsilon), y_{\mathrm{per}}(2 \pi \varepsilon), z_{\mathrm{per}}(2 \pi \varepsilon))$, solution of \eqref{example-2-1} starting at $(0,0,0)$ with controls $u_{1,\mathrm{per}}$ and $u_2$. Integrating system \eqref{example-2-1}, $(x_{\mathrm{per}}(2 \pi \varepsilon), y_{\mathrm{per}}(2 \pi \varepsilon), z_{\mathrm{per}}(2 \pi \varepsilon))$ reads
\begin{equation}
\textstyle \left ( \int_0^{2\pi \varepsilon} u_1^{\mathrm{per}} (t) \mathrm{d} t, \int_0^{2\pi \varepsilon} \!\! \int_0^t u_1^{\mathrm{per}}(\tau) \mathrm{d} \tau \mathrm{d} t , \int_0^{2 \pi \varepsilon} \left ( \left ( \int_0^{t} u_1^{\mathrm{per}}(\tau) \mathrm{d} \tau \right ) ^3 + ( 1 + u_2 (t)) \left ( \int_0^{t}  \int_0^{\tau} u_1^{\mathrm{per}}(\sigma) \mathrm{d} \sigma \mathrm{d} \tau \right ) ^2  \right ) \mathrm{d} t \right ),
\label{eq:xyz}
\end{equation}
which allows, substituting \eqref{eq:u1u2} and \eqref{eq:u1per} in \eqref{eq:xyz}, to explicitly calculate the value of $F(a,b,c)$. With the help of a computer algebra software, we can then calculate the determinant of the Jacobian matrix at $(0,0,0)$ of $F$, which is equal to
\[
-\pi^5 \varepsilon^{14} \left ( \frac{12}{5} \pi^2 \varepsilon - \frac{45}{4} \varepsilon - \frac{184}{27} \pi^2 + \frac{8102}{81} \right ),
\]
which is nonzero for $\varepsilon > 0$. Therefore, $F$ is onto and we apply the inverse mapping theorem: for every $(x,y,z)$ in a neighborhood of $(0,0,0)$, there exists $(a,b,c) \in \mathbb{R}^3$ such that the associated control $(u_{1,\mathrm{per}},u_2)$ drives the system from $(0,0,0)$ to $(x,y,z)$. Hence \eqref{example-2-1} is $\mathrm{W}^{1,\infty}$-STLC at 0. 
\end{xmpl}

This example suggests that the bracket $[f_1,[f_2,f_1]]$ plays a strong role in the controllability of the system, as it allows to recover $\mathrm{W}^{1,\infty}$-STLC that is unattainable in that case with only one scalar control. 

\section{Proofs of the Theorems \label{section:proofs}}

We start with a few notations. Given controls $u_1$, $u_2$ in $\mathrm{L}^{\infty}([0,T])$, we denote by $z_u(T)$ the solution of \eqref{eq:3} with the controls $u_1$, $u_2$ at time $T$, with $z(0)=0$. 

For a vector field $f$ in $\mathcal{X}$ with components $a_1, \dots, a_n$ ($n$  real analytic functions) in coordinates $x=(x_1, \dots, x_n)$, we write indifferently 
\begin{equation*}
f(x) = \begin{pmatrix} a_1(x) \\ \vdots \\ a_n(x) \end{pmatrix} \quad \text{or} \quad f=\sum_{k=1}^n a_k \frac{\partial}{\partial x_k}.
\end{equation*}
Considering $f$ as a differential operator, for a smooth function $\phi$, $f \phi$ is given by $f \phi = \sum_{k=1}^n a_k \frac{\partial \phi}{\partial x_k}$. For $f,g \in \mathcal{X}$, we define the composed operator (of order 2) $fg$ as $(fg)\phi = f(g\phi)$. In coordinates, and if $g = \sum_{k=1}^n b_j \frac{\partial}{\partial x_j}$, one has 
\[
fg=\sum_{k=1}^n \sum_{j=1}^n a_k b_j \frac{\partial^2}{\partial x_k \partial x_j} + \sum_{j=1}^n \left ( \sum_{k=1}^n a_k \frac{\partial b_j}{\partial x_k} \right ) \frac{\partial}{\partial x_j}.
\]
\begin{rmrk}
When $f$ and $g$ are considered as differential operators, their Lie bracket is simply their commutator: $[f,g] = fg-gf$, and it turns out to be a differential operator of order 1 (i.e. a vector field) because the higher order terms cancel.   
\end{rmrk}
For a multi-index $I=(i_1,\ldots,i_k) \in \{0,1,2\}^k$, we denote by $f_I$ the iterated composition of operators $f_{i_1} f_{i_2} \ldots f_{i_k}$ associated to \eqref{eq:3}. 

Let $u=(u_1, u_2)$ in $\mathrm{L}^{\infty}([0,T])$.  For a multi-index $I=(i_1,\ldots,i_k) \in \{0,1,2\}^k$, the iterated integral $\int_0^T u_I $ is defined as
\[
\int_0^T \int_0^{\tau_k} \int_0^{\tau_{k-1}} \dots \int_0^{\tau_2} u_{i_k} (\tau_k) u_{i_{k-1}} (\tau_{k-1} ) \dots u_{i_2} (\tau_2) u_{i_1} (\tau_1) \mathrm{d} \tau_1 \mathrm{d} \tau_2 \dots \mathrm{d} \tau_k,
\]
with the convention $u_0 = 1$.

Let $\Phi : \mathbb{R}^n \rightarrow \mathbb{R}$ a real analytic function defined in a neighbourhood of $0$ in $\mathbb{R}^n$. The Chen-Fliess series associated to $\Phi$ is defined as
\begin{equation}
\Sigma (u,f,\Phi,T) = \sum_{I} \left ( \int_0^T u_I \right ) ( f_I \Phi ) (0).
\label{eq:chen-fliess}
\end{equation}
The summation is made over all the multi-indices $I=(i_1, \dots, i_k)$ in $\{ 0, 1, 2 \}^k$ with $k \in \mathbb{N}^*$. 
The Chen-Fliess series appears in a range of works in control theory and geometry (see \cite{fliess1978developpements,fliess1975series,chen1957integration}). It is shown in \cite[Proposition 4.3, p. 698]{sussmann1983lie} that for all $A>0$, there exists $T_0(A)>0$ such that the series converges for any $T\leqslant T_0$ and $u$ such that $\| u \|_{\mathrm{L}^\infty[0,T]} \leqslant A$, uniformly with respect to $u$ and $T$, to $\Phi(z_u(T))$, i.e. we can write
\[
\Phi(z_u(T)) = \Sigma (u,f,\Phi,T) .
\]
A wisely chosen function $\Phi$ conveniently allows, through the series $\Sigma$, to focus on the brackets of interest, and highlights their role in the following proofs.

\subsection{Proof of Theorem \ref{main-thm-bis}}

The following proof relies on similar arguments than those used to show Proposition \ref{thm:sussmann} in \cite[pp.707-710]{sussmann1983lie}.  We treat both Case \ref{first-case} and Case \ref{second-case} together.

Given $u_2^{\text{eq}}$ in $\mathbb{R}$ verifying moreover $u_2^{\text{eq}} \neq \beta$ in case \ref{first-case}, we can always perform the linear transformation $(u_1,u_2) \mapsto (u_1,u_2 - u_2^{\text{eq}})$ from Remark \ref{remark:S'} to return to the equilibrium $(0,(0,0))$. One can easily check that the transformed system falls within the same case for the equilibrium $(0,(0,0))$ as the original system for the equilibrium $(0,(0,u_2^{\text{eq}}))$. Therefore, we assume from now on that $u_2^{\text{eq}}=0$ and we study controllability around $(0,(0,0))$.

The next step is to define suitable local coordinates, and a real analytic function $\Phi : \mathbb{R}^n \rightarrow \mathbb{R}$ mapping the state $z$ to one of these local coordinates. The function $\Phi$ will then be used to define the system Chen-Fliess series.

We shall prove that we have, for this suitable choice of $\Phi$,
\begin{equation}
\Phi (z_u(T)) \geqslant 0
\label{eq:proof1}
\end{equation}
for any small positive $T$ and any control $u (\cdot ) = (u_1 (\cdot), u_2 (\cdot) )$ satisfying the norm requirements from the STLC definitions (\ref{def-stlc} in case \ref{first-case} or \ref{def-stlcq} in case \ref{second-case}), with $z_u (\cdot)$ defined on the first line of section \ref{section:proofs}. The inequality \eqref{eq:proof1} obviously contradicts local controllability (STLC in case \ref{first-case} or $\alpha$-STLC for any $\alpha$ in case \ref{second-case}). 

Let $d_1=\dim R_1(0)$ and $d_2 = \dim (\mathrm{Span}( [f_1, [f_2, f_1]](0) ) + R_1(0))$ ($d_2$ might be equal to $d_1$ or $d_1+1$). Let $(g_1,\dots,g_n)$ be vector fields in $\mathrm{Lie}(f_0,f_1,f_2)$ such that:
\begin{itemize}
\item $g_1 = f_1$,
\item $(g_1(0),\dots,g_n(0))$ is a basis of $\mathbb{R}^n$,
\item in Case \ref{first-case}, $(g_1(0),\dots,g_{d_1}(0))$ is a basis of $R_1(0)$ and $g_{d_1+1} = [f_1,[f_0,f_1]]$, 
\item in Case \ref{second-case}, $(g_1(0),\dots,g_{d_2}(0))$ is a basis of $\mathrm{Span}( [f_1, [f_2, f_1]](0) ) + R_1(0)$ and $g_{d_2+1} = [f_1,[f_0,f_1]]$.
\end{itemize}
For $s \in \mathbb{R}$ and $g$ a vector field, let $e^{s g}$ denote the flow of $g$ at time $s$. The real analytic map $(s_1, \dots, s_n) \mapsto e^{s_1 g_1} \circ e^{s_2 g_2} \circ \dots \circ e^{s_n g_n} (0)$ sends $0$ to $0$. Moreover, its Jacobian at $0$ is invertible (for its columns are the components of $g_1(0), \dots, g_n(0)$). In virtue of the inverse mapping theorem, it has a real analytic inverse $\zeta \mapsto (s_1 (\zeta), \dots, s_n (\zeta) )$, defined a certain neighborhood $\mathcal{V}$ of $z=0$ in $\mathbb{R}^n$. One then has, for all $\zeta$ in $\mathcal{V}$,
\[
\zeta = e^{s_1(\zeta) g_1} \circ e^{s_2(\zeta) g_2} \circ \dots \circ e^{s_n (\zeta) g_n} (0),
\]
i.e., $(s_1,\dots,s_n)$ are local coordinates for $\zeta$. Then:
\begin{itemize}
\item in case \ref{first-case}, we define $\Phi(\zeta)=s_{d_1+1}(\zeta)$,
\item in case \ref{second-case}, we define $\Phi(\zeta)=s_{d_2+1}(\zeta)$.
\end{itemize}

The function $\Phi$ is real analytical on $\mathcal{V}$ and has, by construction, the following properties:
\begin{align}
\label{prop-phi-1} & \Phi (0)  = 0, \\
\label{prop-phi-2} \text{(Case \ref{second-case})} \quad & \forall g \in R_1, (g \Phi) (0) =0, \\
\label{prop-phi-2bis} \text{(Case \ref{first-case})} \quad & \forall g \in \mathrm{Span}(\{ [f_1, [f_2, f_1]] \}, R_1), (g \Phi) (0) =0, \\
\label{prop-phi-3} & f_1 \Phi = 0 \text{ on } \mathcal{V}, \\
\label{prop-phi-4} &([f_1, [f_1, f_0]] \Phi) (0) = 1. 
\end{align}

We then consider the Chen-Fliess series $\Sigma(u,f,\Phi,T)$ associated to $\Phi$ (see \eqref{eq:chen-fliess}) and split its terms into six different types:
\begin{equation*}
\Sigma (u,f,\Phi,T) = T_1 + T_2 + T_3 + T_4 + T_5 + T_6,
\end{equation*}
where each $T_i$ contains the terms with multi-indices $I$ defined as follows:

\begin{itemize}
\item $T_1 : I=(2,\dots)$, or $I=(0,\dots)$,
\item $T_2 : I=(\dots,1)$,
\item $T_3 : I=(1,J)$ with J containing only 0's and 2's,
\item $T_4 : I=(1,1,0)$,
\item $T_5 : I=(1,1,2)$,
\item $T_6$ :  all the remaining terms.
\end{itemize}

We have $T_1 = T_2 = 0$ because, from \eqref{prop-phi-3} and assumption \eqref{eq:4}, $f_I \Phi (0) = 0$ for $I$ of these types.

We also have $T_3 = 0$; indeed, let $I=(1,i_2, \dots, i_k)$ with $i_j = 0$ or $2$ for all $j \in \{ 2, \dots k \}$. Then we can write that:
\[
f_1 f_{i_2} \dots f_{i_k} = [f_1, f_{i_2} ] f_{i_3} \dots f_{i_k} + f_{i_2} f_1 \dots f_{i_k},
\]
and $(f_{i_2} f_1 \dots f_{i_k} \phi) (0) = 0$ because of assumption \eqref{eq:4} and the fact that $i_2 =0$ or $2$. Similarly, we have
\[
[f_1, f_{i_2} ] f_{i_3} \dots f_{i_k} = [[f_1, f_{i_2}],f_{i_3}] f_{i_4} \dots f_{i_k} + f_{i_3} [f_1, f_{i_2}] \dots f_{i_k},
\]
and $((f_{i_3} [f_1, f_{i_2}] \dots f_{i_k} )\phi ) (0) = 0$ because of assumption \eqref{eq:4} and the fact that $i_3 =0$ or $2$. Repeating this operation $k-3$ more times, we eventually get that 
\[
(f_1 f_{i_2} \dots f_{i_k} \phi) (0) = ([ \dots [f_1, f_{i_2}], \dots, f_{i_k} ] \phi )(0).
\] 
But $[ \dots [f_1, f_{i_2}], \dots, f_{i_k} ]$ is in $R_1$, so $(f_1 f_{i_2} \dots f_{i_k} \phi) (0) =0$ because of \eqref{prop-phi-2}.

In order to calculate $T_4$, we write
\[
\begin{array}{r l}
f_1 f_1 f_0 & = f_1 f_0 f_1 + f_1 [f_1, f_0] \\
 & = f_1 f_0 f_1 - [f_1, f_0] f_1 - [[f_1, f_0], f_1].
 \end{array}
 \]
The first two terms on the right-hand side vanish when evaluated at 0 against $\Phi$ because of assumption \eqref{eq:4}, so $(f_{(1,1,0)} \Phi ) (0) = - ([[f_1, f_0], f_1] \Phi) (0) = 1$ by \eqref{prop-phi-4}.
Moreover, the control integral part is given by 
\[
\begin{array}{r l}
\displaystyle \int_0^T u_{(1,1,0)} & \displaystyle = \int_0^T \int _0^s u_1 (\sigma) \int_0^\sigma u_1 (\tau) \mathrm{d} \tau \mathrm{d} \sigma \mathrm{d} s \\
 & \displaystyle = \int_0^T \int_0^s v_1 ' (\sigma) v_1 (\sigma) \mathrm{d} \sigma \mathrm{d} s \\
 & \displaystyle = \frac{1}{2} \int_0^T v_1^2 (s) \mathrm{d} s 
 \end{array}
 \]
 with 
 \begin{equation}
 v_1(t)=\int_0^t u_1(s) \mathrm{d} s,
 \label{eq:v1}
 \end{equation}
 so overall $T_4 = \frac{1}{2} \| v_1 \|_{\mathrm{L}^2}^2.$
 
For $T_5$, using assumption \eqref{eq:4} again, we obtain
\begin{equation*}
(f_1 f_1 f_2 \Phi) (0) = ( [f_1, [f_1, f_2]] \Phi) (0).
\end{equation*}

This is where the two different cases of the theorem appear. Indeed:
\begin{itemize}
\item Case \ref{second-case}: we have $(f_1 f_1 f_2 \Phi) (0) = 0$ thanks to \eqref{prop-phi-2bis}, so $T_5 = 0$.
\item Case \ref{first-case}:  $[f_1, [f_0, f_1]](0) \in \mathrm{Span}(R_1(0),[f_1, [f_2, f_1]](0))$. 
\end{itemize}
Here, we write that $[f_1, [f_2, f_1]]=- \beta [f_1, [f_0, f_1]] + g$, with $g \in R_1$ and $\beta \in \mathbb{R}^*$. Thanks to \eqref{prop-phi-2} and \eqref{prop-phi-4}, we conclude that
\begin{equation*}
(f_1 f_1 f_2 \Phi) (0) = - \beta. 
\end{equation*}
The control integral associated to $T_5$ reads
\begin{equation*}
\begin{array}{r l}
\int_0^T u_{(1,1,2)}  & = \int_0^T u_2 (s) \int _0^s u_1 (\sigma) \int_0^\sigma u_1 (\tau) \mathrm{d} \tau \mathrm{d} \sigma \mathrm{d} s \\ & \leqslant \frac{1}{2} \, \| u_2 \|_{\mathrm{L}^{\infty}} \| v_1 \|_{\mathrm{L}^2}^2
\end{array}
\end{equation*}
and therefore $T_5$ is bounded:
\begin{equation*}
\left | T_5 \right | \leqslant \frac{1}{2} | \beta | \, \| u_2 \|_{\mathrm{L}^{\infty}} \| v_1 \|_{\mathrm{L}^2}^2.
\end{equation*}

Finally, we are going to show that the terms in $T_6$ add up to a small remainder. Let $I$ be a multi-index such that the associated term in the series is in $T_6$. Then $I = (1, J, 1, K)$ with $K = (k_1, \dots, k_q)$ and $J=(j_1,\dots, j_r)$ such that $q \geqslant 1$, $q + r \geqslant 2$ and $J$ contains only $0$'s and $2$'s. Finally, let us denote by $J_2$ the number of $2$'s in $J$, and $K_1$ and $K_2$ respectively the number of $1$'s and $2$'s in $K$. We write the control integral
\begin{equation}
\int_0^T u_I = \int_0^T \int_0^{s_q} \int_0^{s_{q-1}} \dots \int_0^{s_2} u_{k_q} (s_q) u_{k_{q-1}} (s_{q-1} ) \dots u_{k_1} (s_2) W (s_1) \mathrm{d} s_1 \dots \mathrm{d} s_q,
\label{t6-control}
\end{equation}
with
\begin{equation*}
W (s) = \int_0^s u_1 (\tau_{_r+1}) \int_0^{\tau_{r+1}} u_{j_r} (\tau_r) \dots \int_0^{\tau_1} u_{j_1} (\tau_1) \int_0^{\tau_0} u_1(\tau_0) \mathrm{d} \tau_0 \dots \mathrm{d} \tau_{r+1} .
\end{equation*}
Then, bounding $u_2$ by $\| u_2 \|_{\mathrm{L}^{\infty}}$ and using $u_0 = 1$, we have:
\begin{equation*}
| W (s) | \leqslant \| u_2 \| _{\mathrm{L}^{\infty}}^{J_2} \left | \int_0^s u_1 (\tau_2) \int_0^{\tau_2} \dots \int_0^{\tau_{r+1}} u_1(\tau_{r+2}) \mathrm{d} \tau_{r+2} \dots \mathrm{d} \tau_2 \right |,
\end{equation*}
which reads
\begin{equation}
| W (s) | \leqslant \| u_2 \| _{\mathrm{L}^{\infty}}^{J_2}\left | w_r(s) \right |,
\label{eq:w_r4}
\end{equation}
with 
\begin{equation}
w_r(s) = \frac{1}{r!} \int_0^s \int_0^{\tau} u_1 (\tau) (\tau - \sigma)^r u_1 (\sigma) \mathrm{d} \sigma \mathrm{d} \tau.
\label{eq:w_r5}
\end{equation}
Substituting \eqref{eq:w_r4} in \eqref{t6-control}, and bounding again $u_1$ and $u_2$ respectively by $\| u_1 \|_{\mathrm{L}^{\infty}}$ and $\| u_2 \|_{\mathrm{L}^{\infty}}$, we have:
\begin{equation}
\left | \int_0^T u_I \right | \leqslant \| u_2 \|_{\mathrm{L}^{\infty}}^{K_2+J_2} \| u_1 \|_{\mathrm{L}^{\infty}}^{K_1} \frac{1}{(q-1)!} \int_0^T (T-s)^{q-1} | w_r (s) | \mathrm{d} s.
\label{eq:w_r7}
\end{equation}

The study of $w_r (s)$ splits in three cases:
\begin{itemize}
\item if $r=0$, we have 
\begin{equation}
w_0 (s) = \frac{1}{2} v_1^2 (s),
\label{eq:w_r1}
\end{equation}
where $v_1$ is defined in \eqref{eq:v1}.
\item if $r=1$, we integrate by parts two times \eqref{eq:w_r5} to get 
\begin{equation}
w_1 (s) = v_1 (s) \int_0^s v_1 (\tau) \mathrm{d} \tau - \int_0^s v_1^2 (\tau) \mathrm{d} \tau.
\label{eq:w_r2}
\end{equation}
\item if $r>1$, we integrate by parts two times \eqref{eq:w_r5} to get
\begin{equation}
w_r (s) = - \frac{v_1 (s)}{(r-1)!} \int_0^s (s-\sigma)^{r-1} v_1 (\sigma) \mathrm{d} \sigma + \frac{1}{(r-2)!} \int_0^s v_1 (\tau) \int_0^{\tau} (\tau-\sigma)^{r-2} v_1 (\sigma) \mathrm{d} \sigma \mathrm{d} \tau.
\label{w_r}
\end{equation}
For $r>0$, the Cauchy-Schwarz inequality yields 
\begin{equation}
\left | \int_0^s (s-\sigma)^{r} v_1 (\sigma) \mathrm{d} \sigma \right | \leqslant \frac{s^{r+1/2}}{\sqrt{2r+1}} \| v_1 \|_{\mathrm{L}^2} \quad \text{and} \quad \left | \int_0^s \sigma^{r} v_1 (\sigma) \mathrm{d} \sigma \right | \leqslant \frac{s^{r+1/2}}{\sqrt{2r+1}} \| v_1 \|_{\mathrm{L}^2}.
\label{eq:w_r6}
\end{equation}
The inequalities \eqref{eq:w_r6} applied to \eqref{w_r} lead to the following majoration:
\begin{equation}
| w_r (s) | \leqslant | v_1 (s) | \frac{\| v_1 \|_{\mathrm{L}^2}}{(r-1)!} \frac{s^{r-1/2}}{\sqrt{2r-1}} + \frac{\| v_1 \|_{\mathrm{L}^2}^2}{(r-2)!} \frac{s^{r-1}}{\sqrt{(2r-3)(2r-2)}}.
\label{eq:w_r3}
\end{equation} 
\end{itemize}

Substituting \eqref{eq:w_r1}, \eqref{eq:w_r2}, and \eqref{eq:w_r3} in \eqref{eq:w_r7}, and bounding $(T-s)$ by $T$, we have:
\begin{equation*}
\left | \int_0^T u_I \right | \leqslant \| u_2 \|_{\mathrm{L}^{\infty}}^{K_2+J_2} \| u_1 \|_{\mathrm{L}^{\infty}}^{K_1} \frac{A T^{q+r-1}}{(q-1)!(r-2)!} \| v_1 \|_{\mathrm{L}^2}^2,
\end{equation*}
where $(r-2)!$ is replaced by $1$ if $r \in \{ 0,1 \}$. Here and hereafter, $A$ is a constant that may vary from line to line.

The fields $f_i$ and the function $\Phi$ are real analytic. It is stated in \cite[Lemma 4.2, p.697]{sussmann1983lie} that we have, for some constant $C$ independent of $I$, the majoration
\begin{equation*}
| (f_I \Phi)(0) | \leqslant C^{q+r+2} (q+r+2)!.
\end{equation*}
Hence we can bound the whole term of index $I$ from the series:
\begin{equation*}
\left | \left ( \int_0^T u_I \right ) (f_I \Phi)(0) \right | \leqslant B(q,r) \| T u_2 \|_{\mathrm{L}^{\infty}}^{K_2+J_2} \| u_1 \|_{\mathrm{L}^{\infty}}^{K_1} ,
\end{equation*}
with
\begin{equation*}
B(q,r)=A C^{q+r+2} T^{q+r-2} \frac{(q+r)!(q+r+2)^5}{q!r!}\| v_1 \|_{\mathrm{L}^2}^2,
\end{equation*}
where we bounded $\frac{(q+r+2)}{(q-1)!(r-2)!}$ by $\frac{(q+r)!(q+r+2)^5}{q!r!}$ to encompass the cases $r=0$ and $r=1$.

For any given $q$ and $r$, there are $2^r 3^q$ corresponding indices $I$. More precisely:
\begin{itemize}
\item for any given $q$, $K_2$, $K_1$, there are $\binom{q}{K_1}$ choices to place the $2$'s and then $\binom{q-K_1}{K_2}$ choices to place the $1$'s.  
\item for any given $r$, $J_2$, there are $\binom{r}{J_2}$ choices to place the $1$'s.
\end{itemize}
Therefore, summing all the terms in $T_6$, we obtain a majoration:
\begin{equation}
| T_6 | \leqslant T \sum_{\substack{r \geqslant 0, q \geqslant 1 \\ r+q \geqslant 2}} B(q,r) \sum_{J_2=0}^r \binom{r}{J_2} \| u_2 \|_{\mathrm{L}^{\infty}}^{J_2} \sum_{K_1=0}^{q} \binom{q}{K_1} \| u_1 \|_{\mathrm{L}^{\infty}}^{K_1} \sum_{K_2=0}^{q-K_1} \binom{q-K_1}{K_2} \| u_2 \|_{\mathrm{L}^{\infty}}^{K_2}.
\label{eq:bqr1}
\end{equation}
This rewrites as
\begin{equation}
| T_6 | \leqslant T \sum_{\substack{r \geqslant 0, q \geqslant 1 \\ r+q \geqslant 2}} B(q,r) (1+\| u_2 \|_{\mathrm{L}^{\infty}} )^r (1 + \| u_2 \|_{\mathrm{L}^{\infty}} + \| u_1 \|_{\mathrm{L}^{\infty}} )^q.
\label{eq:bqr3}
\end{equation}
Using $(1+\| u_2 \|_{\mathrm{L}^{\infty}} ) \leqslant (1 + \| u_2 \|_{\mathrm{L}^{\infty}} + \| u_1 \|_{\mathrm{L}^{\infty}} )$ and renumbering the terms of the sum for $p\geqslant 2$ and $0 \leqslant r \leqslant p$ (such that $p=r+q$ in equations \eqref{eq:bqr1} and \eqref{eq:bqr3}), one obtains
\begin{equation}
| T_6 | \leqslant T \sum_{p \geqslant 2} \sum_{r=0}^p B(p-r,r) (1 + \| u_2 \|_{\mathrm{L}^{\infty}} + \| u_1 \|_{\mathrm{L}^{\infty}} )^p.
\label{eq:bqr4}
\end{equation}
There exists $T_0$ such that the series in \eqref{eq:bqr4} converges for all $T$ in $[0,T_0]$.

Including its limit for $T=T_0$ and the other constants in a new constant $D(T_0)$, we finally obtain, for any $T \in [0, T_0]$,
\begin{equation}
| T_6 | \leqslant T D(T_0) \| v_1 \|_{\mathrm{L}^2}^2.
\label{eq:end1}
\end{equation}

We can now end the proof of non-controllability in both cases.

In case \ref{second-case}, we obtain that
\begin{equation}
\Sigma (u,f,\Phi,T) = \frac{1}{2} \| v_1 \|^2_2 + T_6.
\label{eq:end2}
\end{equation}
Let $\varepsilon_0$ be a real positive number such that $\varepsilon D(T_0) \leqslant \frac{1}{2}$. Let $\varepsilon = \min (T_0,\varepsilon_0)$. Using \eqref{eq:end1} in \eqref{eq:end2}, we obtain that $\Sigma (u,f,\Phi,T) \geqslant 0$ for all $T\leqslant \varepsilon$, i.e. we have proven \eqref{eq:proof1}. Hence, system \eqref{eq:3} is not $\alpha$-STLC at $(0,(0,0))$ for any $\alpha \geqslant 0$.

In case \ref{first-case}, we have
\begin{equation}
\Sigma (u,f,\Phi,T) = \frac{1}{2} \| v_1 \|^2_2 + T_5 + T_6,
\label{eq:end3}
\end{equation}
knowing that 
\begin{equation}
|T_5 + T_6| \leqslant \frac{1}{2} | \beta | \| u_2 \|_{\mathrm{L}^{\infty}} \| v_1 \|_{\mathrm{L}^2}^2 + T D(T_0) \|v_1 \|_{\mathrm{L}^2}^2. 
\label{eq:end4}
\end{equation}
Let $\varepsilon_0$ be a real positive number such that $\varepsilon ( \frac{1}{2} | \beta | + D(T_0) ) \leqslant \frac{1}{2}$. Let $\varepsilon = \min (T_0,\varepsilon_0)$. Assume $T \leqslant \varepsilon$ and $\| u_2 \|_{\mathrm{L}^{\infty}} \leqslant \varepsilon$. Using \eqref{eq:end4} in \eqref{eq:end3}, we obtain that $\Sigma (u,f,\Phi,T) \geqslant 0$ for all $T\leqslant \varepsilon$ and $\| u_2 \|_{\mathrm{L}^{\infty}} \leqslant \varepsilon$, i.e. we have proven \eqref{eq:proof1}. Hence, system \eqref{eq:3} is not STLC at $(0,(0,0))$.

This ends the proof of Theorem \ref{main-thm-bis}.

\subsection{Sketch of proof of Theorem \ref{thm:b2}}

The proof is a close adaptation of the proof of Proposition \ref{thm:kawski} given in \cite[pp.40-72]{kawski1986nilpotent} so we only present the main arguments here.

Let $J$ be the set 
\begin{equation*}
J = \{ 0,2 \}^3 \backslash \{ (0,0,0) \}.
\end{equation*}
Like in the proof of Theorem \ref{main-thm-bis}, we build a function $\Phi$ on a neighborhood $\mathcal{V}$ of $0 \in \mathbb{R}^n$  such that:
\begin{align}
\label{prop-phi-1-2} & \Phi (0)=0, \\
\label{prop-phi-2-2} \text{(Case \ref{first-case})} \quad & \forall g \in R_1, (g \Phi) (0) =0, \\
\label{prop-phi-2bis-2} \text{(Case \ref{second-case})} \quad & \forall g \in \mathrm{Span} ( R_1(0), \{ [[f_i,f_1],[f_j,[f_k,f_1]]], (i,j,k) \in J \} ), (g \Phi) (0) =0, \\
\label{prop-phi-3-2} & f_1 \Phi = 0 \text{ on } \mathcal{V}, \\
\label{prop-phi-4-2} & ([f_1, [f_0, f_1]] \Phi) (0) = 1, \\
\label{prop-phi-5-2} & ([f_1, [f_2, f_1]] \Phi) (0) = 0.
\end{align}

Then, we consider the Chen-Fliess series
\[
\sum_{I} \left ( \int_0 ^T u_I \right ) (f_I \Phi) (0),
\]
for some $T>0$ and controls $u_1$, $u_2$ in $\mathrm{L}^{\infty}([0,T])$.

In the previous case, we wanted to show that under certain smallness conditions on the controls, the term associated to $B_1$ (the $T_4$ term) dominates the rest of the series. We were able to do so by, on the one hand, building the function $\Phi$ that makes all the required ``low-order'' terms ($T_1$, $T_2$ and $T_3$) vanish, and, on the other hand, checking that the ``high-order'' terms ($T_6$) add up to a small remainder. The term $T_5$ was the one for which two cases appeared depending on the bracket $B_1$ behavior.

The same path can be followed here: we want to show that the term associated to $B_2$ dominates the rest of the series under certain norm conditions on the controls. 
However, the classification in different types that is used in the proof of Theorem \ref{main-thm-bis} is not as straightforward here. Indeed, higher-order terms in the series are involved in comparison to $B_2$, which greatly increases combinatorial and computational complexity. One can deal with this issue by performing what is called a ``change of basis'' in \cite{kawski1986nilpotent} over the differential operators $f_I$ in the series that are associated with an index $I$ with one, two or three 1's in it, combined with successive integrations by parts of the iterated integrals $\int u_I$. This suitable change of basis is performed in detail in \cite[pp.41-56]{kawski1986nilpotent} for a scalar-input system. 

For $I$ a multi-index, let $k$ be the number of 1's in $I$. The change of basis and integrations by parts eventually allow to replace the part of the series ranging over $I$ with $k \leqslant 3$ by
\begin{equation*}
\sum_{k=1}^3 \sum_{l=0}^{\infty} \sum_{m=0}^l \sum_{J(k,l,m)} b_{J(k,l,m)}^{k,l} ( W_{J(k,l,m)}^{k,l} \Phi ) (0)
\end{equation*}
where $b_{J(k,l,m)}^{k,l}$ is a product of iterated integrals, $W_{j}^{k,l}$ is a differential operator, and $J$ is an index used to enumerate the elements for each triplet of indices $(k,l,m)$ ($k$ refers to the number of 1's, $m$ to the number of 2's and $l$ to the number of 0's and 2's in $J$). More precisely:
\begin{itemize}
\item for $k=1$, $J(1,l,m)$ is a single index $\lambda$ that ranges in $\{ 1 \dots \binom{l}{m} \}$,
\item for $k=2$, $J(2,l,m)$ is a double index $(\lambda, \mu)$ that ranges in $\{ 1, \dots \binom{l}{m} \} \times \{ 1, \dots l+1 \}$,
\item for $k=3$, $J(3,l,m)$ is a quadruple index $(\lambda, \xi, \mu, \nu)$ with $\lambda$ ranging in $\{ 1 \dots \binom{l}{m} \}$, $\xi$ ranging in $\{1,2,3\}$ and $\mu$ and $\nu$ ranging differently depending on $\xi$ and divisibility criteria on $l$; see \cite{kawski1986nilpotent} for more details.
\end{itemize}
Let $l \in \mathbb{N}$. For a given $m \in \{0, \dots , l \}$, each value of $\lambda$ is associated to a multi-index in $\{ 0,2 \}^l$ with exactly $m$ times 2. In particular, for $m=0$, $\lambda$ only takes the value 1 and is associated to $(0,\dots,0)$. The explicit expressions of the differential operators $\mathrm{W}^{k,l}_{J(k,l,0)}$ read
\begin{align}
\label{w1} \mathrm{W}^{1,l}_{1}& = \mathrm{ad}^l_{f_0} f_1, \\
\label{w21} \mathrm{W}^{2,l}_{1,\mu} &= \mathrm{ad}^{l-2 \mu +1}_{f_0} (\mathrm{ad}^{2}_{\mathrm{ad}^{\mu-1}_{f_0} f_1} f_0) \; \text{ if } \;  1 \leqslant 2 \mu -1 \leqslant l, \\
\label{w22} \mathrm{W}^{2,l}_{1,\mu} &= ( \mathrm{ad}^{\mu-1}_{f_0} f_1)(\mathrm{ad}^{l-\mu+1}_{f_0} f_1) \; \text{ if } \; l/2 < \mu < l+1, \\
\label{w31} \mathrm{W}^{3,l}_{1,1,\nu,\mu} &= - (\mathrm{ad}^{l-1-\mu-2\nu}_{f_0} [\mathrm{ad}^{\mu}_{f_0} f_1, (\mathrm{ad}^{2}_{\mathrm{ad}^{\nu}_{f_0} f_1} f_0)]) \; \text{ if } \; 0 \leqslant \nu \leqslant \mu \; \text{ and } \; 2\nu + \mu \leqslant l - 1, \\
\label{w32} \mathrm{W}^{3,l}_{1,2,\nu,\mu} &= (\mathrm{ad}^{l-1-\nu-2\mu}_{f_0} (\mathrm{ad}^{2}_{\mathrm{ad}^{\mu}_{f_0} f_1} f_0))(\mathrm{ad}^{\nu}_{f_0} f_1) \; \text{ if } \; 0 \leqslant \nu, \mu \; \text{ and } \; 2\mu + \nu \leqslant l-1, \\
\label{w33} \mathrm{W}^{3,l}_{1,3,\nu,\mu} &= (\mathrm{ad}^{l-\nu-\mu}_{f_0} f_1)(\mathrm{ad}^{\nu}_{f_0} f_1)(\mathrm{ad}^{\nu}_{f_0} f_1) \; \text{ if } \; 0 \leqslant \nu \leqslant \mu \leqslant l-\nu - \mu.
\end{align}
For $m > 0$, the other operators $\mathrm{W}^{k,l}_{J(k,l,m)}$ are obtained by replacing the set of $l$ zeros in the expressions \eqref{w1} to \eqref{w33} by every possible multi-index in $\{ 0,2 \}^l$ containing exactly $m$ times 2. 

The generic expressions of the iterated integrals $b_{j}^{k,l}$ are available in \cite{kawski1986nilpotent}. 

Observe that, by the properties of $\Phi$:
\begin{itemize}
\item for $l \geqslant 0$ and for all $\lambda$, $(\mathrm{W}^{1,l}_{1} \Phi)(0) = 0$,
\item for $l \geqslant 1$ and for all $\lambda$, $(\mathrm{W}^{2,l}_{\lambda,1} \Phi)(0) = 0$,
\item for $l \geqslant 2$ and for all $\lambda$, $(\mathrm{W}^{2,l}_{\lambda,l} \Phi)(0) = 0$,
\item for $l \geqslant 0$ and for all $\lambda$, $(\mathrm{W}^{2,l+1}_{\lambda,l+1} \Phi)(0) = 0$,
\item for $\mu \geqslant 0$, $l \geqslant 2 \mu$, and for all $\lambda$, $(\mathrm{W}^{3,l}_{\lambda,2,\mu, 0} \Phi)(0) = 0$,
\item for $l \geqslant 2$ and for all $\lambda$, $(\mathrm{W}^{3,l}_{\lambda,2,0,1} \Phi)(0) = 0$,
\item for $\mu \geqslant 0$, $l \geqslant 2 \mu$, and for all $\lambda$, $(\mathrm{W}^{3,l}_{\lambda,3,\mu,0} \Phi)(0) = 0$,
\item for $l \geqslant 3$, and for all $\lambda$, $(\mathrm{W}^{3,l}_{\lambda,3,1,1} \Phi)(0) = 0$.
\end{itemize}

Moreover, the term associated to $B_2$ is $\mathrm{W}^{2,3}_{1,2}$, which reads $$\frac{1}{2} \int_0^T \left ( \int_0^t \int_0^s u_1 (r) \mathrm{d} r \mathrm{d} s \right ) ^2 \mathrm{d} t,$$ that we will write $\frac{1}{2} \int (\iint u_1)^2$ to lighten notations.

The next step is to examine the terms associated to the seven brackets $[[f_i,f_1],[f_j,[f_k,f_1]]]$ with $(i,j,k) \in J$, i.e., to the $\mathrm{W}^{2,3}_{\lambda, 2}$ for all values of $\lambda$ and for $m=1,2,3$.  For $(i,j,k)$ in $J$, let 
\begin{equation}
F_{ijk}= [[f_i,f_1],[f_j,[f_k,f_1]]].\end{equation}

Case \ref{prop-first-case} : $B_2(0) \in \mathrm{Span} ( R_1(0), \{ F_{ijk}, (i,j,k) \in J \} )$. This means that there exists $(i,j,k) \in J $ such that $F_{ijk} \Phi (0) \neq 0$. Let $\beta = \max_{J} \big | F_{ijk} \Phi (0) \big | $. The control integrals associated to the $F_{ijk}$ are respectively given by 
\begin{equation*} \begin{array}{c}
\textstyle \iint (\int u_1) (\int u_2 \int u_1), \quad \iint u_2 (\int u_1) (\iint u_1 ), \quad \frac{1}{2} \int u_2 (\iint u_1)^2, \\ \frac{1}{2} \int (\int u_2 \int u_1)^2, \quad \int u_2 \int (\int u_1) (\int u_2 \int u_1), \quad \int u_2 \int u_2 (\int u_1) (\iint u_1), \\ \frac{1}{2} \int u_2 (\int u_2 \int u_1)^2.
\end{array}
\end{equation*}
The first three are bounded in absolute value by $\frac{1}{2}  \| u_2 \|_{\mathrm{L}^{\infty}} \int (\iint u_1)^2$, the following three are bounded in absolute value by 
 $\frac{1}{2}  \| u_2 \|_{\mathrm{L}^{\infty}}^2 \int (\iint u_1)^2$, and the last one is bounded in absolute value by $\frac{1}{2}  \| u_2 \|_{\mathrm{L}^{\infty}}^3 \int (\iint u_1)^2$. The sum of the seven corresponding terms in the series, called $\tau$, is then bounded accordingly:
 \[
 | \tau | \leqslant \frac{1}{2} \beta \| u_2 \|_{\mathrm{L}^{\infty}} (3 + 3  \| u_2 \|_{\mathrm{L}^{\infty}} +  \| u_2 \|_{\mathrm{L}^{\infty}}^2) \textstyle \int (\iint u_1)^2.
 \]
 If we assume that  $\| u_2 \|_{\mathrm{L}^{\infty}}$ is small enough, then $| \tau |$ is negligible (i.e. less than $C \int (\iint u_1)^2$ for any $C>0$).

Case \ref{prop-second-case} : if $B_2(0)$ is not in $\mathrm{Span} ( R_1(0), \{ F_{ijk}, (i,j,k) \in J \} )$, then we observe that for every $(i,j,k)$ in $J$, $F_{ijk} \Phi (0) = 0$ thanks to \eqref{prop-phi-2bis-2}.

For both cases, it is shown in \cite[Proposition 6]{beauchard2018quadratic}, using the Gagliardo-Nirenberg inequality, that 
\[
\| u_1 \|^3_{\mathrm{L}^3} \leqslant c \iint \| u_1 \|_{\mathrm{L}^2}^2 \| u_1 \|_{\mathrm{W}^{1,\infty}},
\]
for a constant $c$ independent of $T$. This is used to show that the sum of the terms associated to $\mathrm{W}^{3,l}_{\lambda,1,0,0}$ for $l\geqslant 1$ and for all $\lambda$ is negligible compared to $\frac{1}{2} \int (\iint u_1)^2$ (i.e. less than $C \int (\iint u_1)^2$ for any $C>0$ in absolute value if $T$ is sufficiently small), provided $\| u_1 \|_{\mathrm{W}^{1,\infty} ([0,T])} \leqslant \varepsilon$. 

It is shown in detail in \cite[pp.62-72]{kawski1986nilpotent} that the remaining terms for $k=2$ and $k=3$ are negligible compared to $\frac{1}{2} \int (\iint u_1)^2$ for a scalar-input system. 

Finally, it is shown in detail in \cite[proof of Property (P)]{stefani1986local} that the ``high-order" remaining terms in the series, associated with indices containing more that four 1's, are negligible compared to $\frac{1}{2} \int (\iint u_1)^2$ for a scalar-input system. 

The computations found in those are straightforwardly adaptable to the two-control case, by systematically majoring $u_2(t)$ with $\| u_2 \|_{\mathrm{L}^{\infty}}$ in the iterated integrals.

We conclude as in the proof of Theorem \ref{main-thm-bis} : gathering all the terms, we see that $B_2$ dominates the series, which means that $\Phi(z_u(T)) \geqslant 0$, and therefore system \eqref{eq:3} is not controllable under the norm conditions on the controls specified in cases \ref{prop-first-case} and \ref{prop-second-case}.
 
\begin{rmrk}
As stated above, the result from Theorem \ref{thm:b2} is only partial, for it does not deal with the case where $B_2(0) \in \mathrm{Span}([f_1,[f_2,f_1]](0),R_1(0))$. However, one cannot ensure that the term associated to the bracket $[f_1,[f_2,f_1]]$ in the Chen-Fliess series is dominated by the $B_2$ term, no matter which smallness hypothesis are made on the controls. The method used in our proofs is therefore not applicable in that particular case. 

Example \ref{example:2-3} highlights the role of the $[f_1,[f_2,f_1]]$ bracket, for it displays a system that is indeed $\mathrm{W}^{1,\infty}$-STLC thanks to $[f_1,[f_2,f_1]]$. By analogy with Theorem \ref{main-thm-bis} for the first obstruction, one could conjecture that the system is not $\mathrm{W}^{2,\infty}$ in this case. Nonetheless, we have not found yet any example that could confirm or infirm this hypothesis.
\end{rmrk}

\section{Conclusion \label{section:conclusion}}

In this paper, we explored the controllability properties of systems with two controls \eqref{eq:3}, satisfying assumption \eqref{eq:4}. 
We have stated two results on these systems, Theorem \ref{main-thm-bis} and Theorem \ref{thm:b2}, providing necessary conditions for local controllability around equilibria. 
These results extend the classical necessary conditions stated for scalar-control systems in \cite{sussmann1983lie}. 
Moreover, they are, to the best of our knowledge, the first results of this nature for non-scalar-input systems. 

This work does not only present a theoretical interest for control theory. 
Using Theorem \ref{main-thm-bis}, we were able to solve the open question of the local controllability of magnetically controlled micro-swimming robots (see paragraph \ref{sssection:microswimmer}). 
One can use our results to easily and systematically address local controllability issues in similar applied situations. 

Our necessary conditions are only based on the ``bad'' brackets $B_1$ and $B_2$ (see equations \eqref{eq:first-bad} and \eqref{eq:second-bad}), but there are higher-order brackets that may prevent $S_2(0)$ to be contained in $S_1(0)$ (in the single-input case see for instance \cite{beauchard2018quadratic}). 
Giving necessary conditions based on these brackets, for instance adapting the results from \cite{beauchard2018quadratic} to the situation \eqref{eq:3}-\eqref{eq:4}, is a possible continuation of the present work. The complexity of the higher-order terms structure in the Chen-Fliess series however makes the analysis very intricate.

\bibliography{biblio}
\bibliographystyle{IEEEtran}

\end{document}